\documentclass[12pt]{article}
\usepackage{amsmath}
\usepackage{amssymb}
\usepackage[dvips]{graphicx}
\usepackage{epsfig}
\usepackage{vector}
\def\eps{\varepsilon}
\font\tencmmib=cmmib10 \skewchar\tencmmib '60
\newfam\cmmibfam
\textfont\cmmibfam=\tencmmib

\def\bbox{\quad\hbox{\vrule \vbox{\hrule \vskip2pt \hbox{\hskip2pt
\vbox{\hsize=1pt}\hskip2pt} \vskip2pt\hrule}\vrule}}
\def\lessim{\ \lower4pt\hbox{$
\buildrel{\displaystyle <}\over\sim$}\ }
\def\gessim{\ \lower4pt\hbox{$\buildrel{\displaystyle >}
\over\sim$}\ }

\def\E{{\cal E}}

\def\M{{\cal M}}
\def\G{{\Gamma}}
\def\eps{{\varepsilon}}
\def\ch{{\mbox{ch}}}

\def\th{{\mbox{th}}}
\def\Av{{\mbox{\rm{Av}}}}

\def\la{{\Bigl\langle}}
\def\ra{{\Bigr\rangle}}

\def\qed{\hfill\break\rightline{$\bbox$}}
\newcommand{\e}{\mathbb{E}}
\newcommand{\p}{\mathbb{P}}
\newcommand{\Reals}{\mathbb{R}}
\newcommand{\Natural}{\mathbb{N}}

\newcommand{\F}{{\cal F}}
\newcommand{\vsi}{{\vec{\sigma}}}

\newtheorem{proposition}{Proposition}

\newtheorem{theorem}{Theorem}

\makeatletter
\@addtoreset{equation}{section}

\makeatother

\font\tencmmib=cmmib10 \skewchar\tencmmib '60
\newfam\cmmibfam
\textfont\cmmibfam=\tencmmib

\def\bbox{\quad\hbox{\vrule \vbox{\hrule \vskip2pt \hbox{\hskip2pt
\vbox{\hsize=1pt}\hskip2pt} \vskip2pt\hrule}\vrule}}
\def\lessim{\ \lower4pt\hbox{$
\buildrel{\displaystyle <}\over\sim$}\ }
\def\gessim{\ \lower4pt\hbox{$\buildrel{\displaystyle >}
\over\sim$}\ }

\def\eps{\varepsilon}

\def\go0{\to 0}

\def\la{\langle}

\def\leftitem#1{\item{\hbox to\parindent{\enspace#1\hfill}}}

\def\qed{\hfill\break\rightline{$\bbox$}}

\def\ra{\rangle}

\def\sg{\sigma}

\def\sg2{\sigma^2}

\def\__{_{\infty}}
\begin{document}

\title{
Bounds for diluted mean-fields spin glass models
}

\author{ 
Dmitry Panchenko\thanks{Department of Mathematics, Massachusetts Institute
of Technology, 77 Massachusetts Ave, Cambridge, MA 02139
email: panchenk@math.mit.edu}, 
Michel Talagrand\thanks{
Equipe d'Analyse de l'Institut Math\'ematique, 4 place Jussieu, 
75230 PARIS Cedex 05 and Department of Mathematics, 
the Ohio-State university, Columbus, OH 43210. 
email: talagran@math.ohio-state.edu\hfill
Supported by NSF grant}
}

\maketitle

\begin{abstract}
  In an important recent paper, \cite{FL}, S. Franz and M. Leone prove rigorous lower bounds for the free energy of the diluted $p$-spin model and the $K$-sat model at any temperature. We show that the  results for these two models are consequences of a single general principle. Our calculations are significantly simpler than those of \cite{FL}, even in the replica-symmetric case. 
\end{abstract}
\vspace{0.5cm}


\section{Introduction.}

Let $p\geq 2$ be an even integer that will be fixed throughout this paper.
For $N\geq 1$ let $\Sigma_N=\{-1,+1\}^N.$
Consider a random function $\theta: \{-1,+1\}^p \to \Reals$ and 
a sequence $(\theta_k)_{k\geq 1}$ of independent copies of $\theta.$
Consider an i.i.d. sequence of indices $(i_{l,k})_{l,k\geq 1}$ with uniform 
distribution on $\{1,\ldots,N\},$ and let 
$M$ be a Poisson r.v. with mean $\e M = \alpha N.$
Let us define the Hamiltonian $H_N(\vsi)$ on 
$\Sigma_N$ by
\begin{equation} 
-H_N(\vsi)=\sum_{k\leq M}\theta_k(\sigma_{i_{1,k}},
\ldots,\sigma_{i_{p,k}})
+H'_N(\vsi),
\label{Ham}
\end{equation}
where $H_N'(\vsi)$ is an arbitrary random function on $\Sigma_N$
independent of all other r.v. in (\ref{Ham}).
The main goal of this paper is to prove upper bounds for
$$
F_N=\frac{1}{N}\e \log \sum_{\vsi\in\Sigma_N} \exp(-H_N(\vsi)).
$$
We will make the following assumptions on the random function $\theta.$
We assume that there exists a random function $f:\{-1,+1\}\to \Reals$ 
such that
\begin{equation}
\exp \theta(\sigma_1,\ldots,\sigma_p) = a(1+ b f_1(\sigma_1)\ldots
f_p(\sigma_p)),
\label{condition1}
\end{equation}
where $f_1,\ldots,f_p$ are independent copies of $f,$
$b$ is a r.v. independent of $f_1,\ldots,f_p$ that satisfies the condition
\begin{equation}
\forall n\geq 1\,\,\,\,\,\,\, \e (-b)^n \geq 0,
\label{condition2}
\end{equation}
and $a$ is an arbitrary r.v.
Finally, we assume that
\begin{equation}
|bf_1(\sigma_1)\ldots f_p(\sigma_p)|< 1 \mbox{ a.s. }
\label{condition3}
\end{equation}

Let us consider two examples when the conditions
(\ref{condition1}), (\ref{condition2}) and (\ref{condition3}) are satisfied. 

{\bf Example 1.} ($p$-spin model) Consider $\beta>0$ and 
a symmetric r.v. $J.$ The $p$-spin model corresponds to
the choice of
$$
\theta(\sigma_1,\ldots,\sigma_p)= \beta J \sigma_1\ldots \sigma_p.
$$ 
(\ref{condition1}) holds with
$a=\ch (\beta J),$ $b = \th (\beta J)$ and $f(\sigma)=\sigma$
and condition (\ref{condition2}) holds
since we assume that the distribution of $J$ is symmetric.

{\bf Example 2.} ($K$-sat model) 
Consider $\beta>0$ and 
a sequence of i.i.d. Bernoulli r.v. $(J_l)_{l\geq 1}$
with $\p(J_l=\pm 1)=1/2.$
The $K$-sat model corresponds to 
$$
\theta(\sigma_1,\ldots,\sigma_p)=-\beta\prod_{l\leq p} \frac{1+J_l \sigma_l}{2}.
$$
(\ref{condition1}) holds with
$a=1,$ $b=e^{-\beta}-1$ and $f_l(\sigma_l)=(1+J_l\sigma_l)/2$
and (\ref{condition2}) holds since $b<0.$

We now introduce certain quantities that will play a fundamental role in the paper. Given a function $f:\Reals^p\to\Reals$ and a vector
$\vec{x}=(x_1,\ldots,x_p)$ let us define
$$
\la f \ra_{\vec{x}}^{-}=
\frac{\sum_{\eps_1,\ldots,\eps_{p-1}=\pm 1} 
f(\eps_1,\ldots,\eps_p)
\exp\sum_{l=1}^{p-1} x_l \eps_l}
{\sum_{\eps_1,\ldots,\eps_{p-1}=\pm 1} 
\exp\sum_{l=1}^{p-1} x_l \eps_l}
$$
(so that $\la f \ra_{\vec{x}}^{-}$ implicitly depends on the last coordinate
$\eps_p$)
and
$$
\la f \ra_{\vec{x}}=
\frac{\sum_{\eps_1,\ldots,\eps_p=\pm 1} 
f(\eps_1,\ldots,\eps_p)
\exp \sum_{l=1}^{p} x_l \eps_l}
{\sum_{\eps_1,\ldots,\eps_{p}=\pm 1} 
\exp\sum_{l=1}^{p} x_l \eps_l}.
$$
Let us define
$$
\E(\eps_1,\ldots,\eps_{p})=\exp \theta (\eps_1,\ldots,\eps_{p}).
$$ 
If the condition (\ref{condition1}) holds then 
\begin{eqnarray}
&&
\la\E\ra_{\vec{x}}^{-}=
\la \exp \theta \ra_{\vec{x}}^{-} 
=
\Bigl\langle a\bigl(1+bf_1(\eps_1)\ldots f_p(\eps_p)\bigr) 
\Bigr\rangle_{\vec{x}}^{-}
\label{E1}
\\
&&
=
a\Bigl(1+b f_p(\eps_p)\prod_{l=1}^{p-1} 
\frac{\sum_{\eps_l=\pm 1} f_l(\eps_l)\exp \eps_l x_l}
{\sum_{\eps_l=\pm 1} \exp \eps_l x_l}\Bigr)
=
a\Bigl(1+b f_p(\eps_p)\prod_{l=1}^{p-1} 
\frac{\Av f_l(\eps)\exp \eps x_l}
{\ch (x_l)}\Bigr),
\nonumber
\end{eqnarray}
where $\Av$ means average over $\eps=\pm 1$
and, similarly,
\begin{equation}
\la \E \ra_{\vec{x}} =
a\Bigl(1+ b \prod_{l\leq p}\frac{\Av f_l(\eps)\exp (x_l \eps)}{\ch(x_l)}\Bigr).
\label{E2}
\end{equation}
Finally, let us define
\begin{equation}
U(\theta,x_1,\ldots,x_{p-1},\eps)=\log
\la \E \ra_{\vec{x}}^{-}\bigr|_{\eps_p=\eps}.
\label{U}
\end{equation}
In the case of the $p$-spin model,  we have
$$
\la \E \ra_{\vec{x}}^{-} =
\ch(\beta J)\Bigl(1+\th(\beta J)\eps_p \prod_{l\leq p-1}\th(\beta x_l)\Bigr)
$$
and
$$
\la \E \ra_{\vec{x}} =
\ch(\beta J)\Bigl(1+\th(\beta J) \prod_{l\leq p}\th(\beta x_l)\Bigr).
$$
In the case of the $K$-sat model, we have
$$
\la \E \ra_{\vec{x}}^{-} =
1+(e^{-\beta}-1)\frac{1+J_p \eps_p}{2}
\prod_{l\leq p-1}\frac{1+J_l\th(\beta x_l)}{2}
$$
and
$$
\la \E \ra_{\vec{x}} =
1+(e^{-\beta}-1)
\prod_{l\leq p}\frac{1+J_l\th(\beta x_l)}{2}.
$$

\section{The replica-symmetric bound.}

Given an arbitrary probability measure $\zeta$ on $\Reals$
consider an i.i.d. sequence $x_l^{i,j},$   $i,j,l\geq 1$ with
distribution $\zeta$ and consider
$U_{i,j}(\eps)=U(\theta_{i,j},x_1^{i,j},\ldots,x_{p-1}^{i,j},\eps),$
where $\theta_{i,j}$ are independent copies of $\theta.$
Let us consider the Hamiltonian
\begin{equation} 
-H_{N,t}(\vsi)=\sum_{k\leq M_t}\theta_k(\sigma_{i_{1,k}},
\ldots,\sigma_{i_{p,k}}) + 
\sum_{i\leq N}\sum_{j\leq k_{i,t}} U_{i,j}(\sigma_i)
+H'_N(\vsi),
\label{inter-rs}
\end{equation}
where $M_t$ is a Poisson r.v. with mean $\e M_t = t \alpha N,$ 
and $k_{i,t}, i\leq N$ are i.i.d. Poisson r.v. with mean 
$\e k_{i,t} = (1-t) \alpha p.$
Let us define
$$
\varphi(t)=\frac{1}{N}\e \log \sum_{\vsi\in\Sigma_N} \exp(-H_{N,t}(\vsi)).
$$
Clearly, $F_N = \varphi(1).$ The following Theorem holds.

\begin{theorem}\label{RS} (RS bound).
If conditions (\ref{condition1}),
(\ref{condition2}) and (\ref{condition3}) hold then
\begin{equation}
F_N \leq \Phi(\zeta)=
\varphi(0) - \alpha(p-1) \e\log \la \E \ra_{\vec{x}},
\label{RSbound}
\end{equation}
where $\vec{x}=(x_1,\ldots,x_p)$ is a vector of i.i.d. r.v. with
the distribution $\zeta.$
\end{theorem}

First of all, since $F_N$ does not depend on $\zeta,$
Theorem \ref{RS} implies
\begin{equation}
F_N\leq \Phi_0 = \inf_{\zeta} \Phi(\zeta).
\label{rsb00}
\end{equation}
Even though Theorem \ref{RS} holds for any $H'_N(\vsi),$
it is particularly interesting when
\begin{equation}
H'_N(\vsi)=\sum_{i\leq N} h_i \sigma_i,
\label{H2}
\end{equation}
where $(h_i)_{i\leq N}$ is a sequence of i.i.d. random variables. With this choice we can write
\begin{eqnarray*}
\varphi(0)&=&\frac{1}{N}\e \log \sum_{\vsi\in\Sigma_N} 
\exp\sum_{i\leq N}\Bigl(
\sum_{j\leq k_{i,1}} U_{i,j}(\sigma_i)+h_i\sigma_i\Bigr)
\\
&=&\frac{1}{N}\e \log \prod_{i \leq N} \sum_{\sigma_i = \pm 1}
\exp\Bigl(
\sum_{j\leq k_{i,1}} U_{i,j}(\sigma_i)+h_i\sigma_i\Bigr)
\\
&=&
\log 2 + \e \log \Av\exp\Bigl(\sum_{j\leq k} U_{j}(\eps)+h\eps\Bigr),
\end{eqnarray*}
where $\Av$ mean average over $\eps=\pm 1,$
and $(U_j),$  $h$ and $k$ are copies of $(U_{1,j}),$ $h_1$ and $k_{1,1}$
correspondingly.
In this case
the bound (\ref{RSbound}) is usually  written in Physics in terms of the functions
$B(x_1,\ldots,x_{p-1})$ and $u(x_1,\ldots,x_{p-1})$ defined by
\begin{equation}
\forall \eps_p=\pm 1 \,\,\,\,\,
Be^{\eps_p u}=\la \E \ra_{\vec{x}}^{-}.
\label{Bu}
\end{equation}
In the case of the $p$-spin model (\ref{Bu}) defines
$$
u=\th^{-1}\Bigl( 1+\th(\beta J)\prod_{l\leq p-1} \th(\beta x_{l}) \Bigr)
$$
and $B=\ch(\beta J)/\ch(\beta u).$ In the case of the $K$-sat model
$$
\th u = \Bigl(\frac{b}{2}J_p\prod_{l\leq p-1}\frac{1+J_l \th(\beta x_l)}{2}\Bigr)
\Big/\Bigl(1+\frac{b}{2}\prod_{l\leq p-1}\frac{1+J_l \th(\beta x_l)}{2}\Bigr)
$$
and
$$
B=\frac{1}{\ch u}\Bigl(1+
\frac{b}{2}\prod_{l\leq p-1}\frac{1+J_l \th(\beta x_l)}{2}
\Bigr).
$$
To write  $\varphi(0)$ in terms of these functions $B$ and $u$, we observe that 
\begin{eqnarray}
\varphi(0)&=&
\log 2 + \e \log \Av\exp\Bigl(\sum_{j\leq k} U_{j}(\eps)+h\eps\Bigr)
\nonumber
\\
&=&
\log 2 + \e\log\prod_{j\leq k}B_j 
\Av\exp\Bigl(\sum_{j\leq k} u_{j}\eps+h\eps\Bigr),
\nonumber
\\
&=&
\log 2 + \alpha p \e \log B + \e \log \ch(\sum_{j\leq k} u_{j} +h),
\label{RSBu}
\end{eqnarray}
using that $\e k = \alpha p$ in the last line.

In the case when $H'_N(\vsi)=0,$ 
it was proved in \cite{T1} (see Chapter 7, \cite{SG}) that for $\alpha$ small enough, 
$$\lim_{N\to\infty} F_N = \Phi(\zeta_{\alpha})=\inf_{\zeta}\Phi(\zeta)
$$
where
$\zeta_{\alpha}$ is the unique solution of the equation
$$
x\sim \sum_{j\leq k} u_j,
$$
where $u_j=u_j(x_1^j,\ldots,x_{p-1}^j)$ is defined in (\ref{Bu}),
$x$ and $x_l^j$ are i.i.d. with the distribution $\zeta_{\alpha},$ $k$ is Poisson with $\e k =\alpha p$ and $\sim$ means equality in distribution.

{\bf Proof of Theorem \ref{RS}.}
Let us consider the partition function
$$
Z=\sum_{\vsi\in\Sigma_N} \exp(-H_{N,t}(\vsi))
$$
(for simplicity of notations we omit the dependence of $Z$ on $N$ and $t$)
and define 
$$
Z_{m}=Z\bigr|_{M_t = m}\,\,\,  \mbox{ and } \,\,\,
Z_{i,k}=Z\bigr|_{k_{i,t} = k}.
$$
If we denote the Poisson p.f. as
$\pi(\lambda,k)=(\lambda^k/k!) e^{-\lambda}$
then
$$\e \log Z = \sum_{m\geq 0} \pi(t \alpha N,m) \e\log Z_m$$ 
and, for any $i \leq N$, 
$$\e \log Z = \sum_{k\geq 0} \pi((1-t)\alpha p, k) \e \log Z_{i,k}.$$
Using the notation $I(m\geq1) =1$ if $m\geq 1$ and $I(m\geq1)=0$ if $m=0$, we have 
\begin{eqnarray}
\varphi'(t)
&=&
\sum_{m=0}^{\infty} \frac{\partial\pi(t\alpha N,m)}{\partial t}
\frac{1}{N}\e \log Z_m +
\sum_{i=1}^N \sum_{k=0}^{\infty} 
\frac{\partial\pi((1-t)\alpha p, k)}{\partial t}
\frac{1}{N}\e\log Z_{i,k}
\nonumber
\\
&
=
&
\alpha \sum_{m=0}^{\infty} 
\bigl(\pi(t\alpha N,m-1)I(m\geq 1) - \pi(t\alpha N,m) \bigr)
\e \log Z_m 
\nonumber
\\
&&
- \alpha p
\frac{1}{N}
\sum_{i=1}^N \sum_{k=0}^{\infty} 
\bigl(\pi((1-t)\alpha p, k-1)I(k\geq 1) - \pi((1-t)\alpha p,k)\bigr)
\e\log Z_{i,k}
\nonumber
\\
&
=
&
\alpha \Bigl(\sum_{m=0}^{\infty} 
\pi(t\alpha N,m) \e \log Z_{m+1}
- \e \log Z \Bigr)
\nonumber
\\
&&
- \alpha p
\frac{1}{N}
\sum_{i=1}^N 
\Bigl(
\sum_{k=0}^{\infty} 
\pi((1-t)\alpha p, k) \e\log Z_{i,k+1}
- \e\log Z \Bigr).
\label{dersteps}
\end{eqnarray}
If we denote by $\la\cdot\ra_m$ the averaging w.r.t. the Gibbs
measure corresponding to the Hamiltonian (\ref{inter-rs})
for a fixed $M_t = m$ then
$$
Z_{m+1}=Z_m \Bigl\langle \exp \theta_{m+1}(\sigma_{i_{1,m+1}},\ldots,
\sigma_{i_{p,m+1}})\Bigr\rangle_m
$$
and, therefore,
\begin{eqnarray*}
&&
\sum_{m=0}^{\infty} 
\pi(t\alpha N,m) \e \log Z_{m+1} = 
\sum_{m=0}^{\infty} 
\pi(t\alpha N,m) \e \log Z_{m} 
\\
&&
+
\sum_{m=0}^{\infty} 
\pi(t\alpha N,m) \e \log \Bigl\langle \theta_{m+1}(\sigma_{i_{1,m+1}},\ldots,
\sigma_{i_{p,m+1}}) \Bigr\rangle_m 
\\
&&
=\e \log Z + \frac{1}{N^p}
\sum_{i_1,\ldots,i_p=1}^N
\e \log 
\Bigl\langle \exp \theta(\sigma_{i_{1}},\ldots,\sigma_{i_{p}})\Bigr\rangle,
\end{eqnarray*}
where $\la\cdot \ra$ now denotes
averaging w.r.t. the Gibbs
measure corresponding to the Hamiltonian (\ref{inter-rs})
and $\theta$ is independent of the randomness in $\la\cdot\ra.$
Similarly,
$$
\sum_{k=0}^{\infty} 
\pi((1-t)\alpha p, k) \e\log Z_{i,k+1}
=
\e \log Z +
\e \log \Bigl\langle \exp U(\sigma_i) \Bigr\rangle,
$$
where $U(\sigma_i)=U(x_1,\ldots,x_{p-1},\sigma_i)$ and where
$x_1,\ldots,x_{p-1}$ are independent of randomness in $\la\cdot\ra.$
Finally, (\ref{dersteps}) implies
\begin{equation}
\varphi'(t)=\alpha\Bigl(
\frac{1}{N^p}
\sum_{i_1,\ldots,i_p=1}^N
\e \log 
\Bigl\langle \exp \theta(\sigma_{i_{1}},\ldots,\sigma_{i_{p}})\Bigr\rangle
-p\frac{1}{N}\sum_{i=1}^N
\e \log \Bigl\langle \exp U(\sigma_i) \Bigr\rangle
\Bigr).
\label{deriv}
\end{equation}
Since $\varphi(1)=\varphi(0)+\int_{0}^1 \varphi'(t)dt$ and since
$\e \log \la \E \ra_{\vec{x}}$ does not depend on $t$,
to prove Theorem \ref{RS} it suffices to show that
\begin{equation}
\frac{1}{N^p}
\sum_{i_1,\ldots,i_p=1}^N
\e \log 
\Bigl\langle \exp \theta(\sigma_{i_{1}},\ldots,\sigma_{i_{p}})\Bigr\rangle
-p\frac{1}{N}\sum_{i=1}^N
\e \log \Bigl\langle \exp U(\sigma_i) \Bigr\rangle
+(p-1)\e \log \la \E \ra_{\vec{x}} \leq 0.
\label{neg}
\end{equation}
By assumptions (\ref{condition1}) and (\ref{condition2}) we can write
\begin{eqnarray*}
\log \Bigl\langle \exp \theta(\sigma_{i_{1}},\ldots,\sigma_{i_{p}})\Bigr\rangle
&=&
\log a + 
\log \Bigl(1+b\Bigl\langle f_1(\sigma_{i_1})\ldots f_p(\sigma_{i_p})\Bigr\rangle\Bigr)
\\
&=& \log a - 
\sum_{n=1}^{\infty}\frac{(-b)^n}{n}\Bigl\langle f_1(\sigma_{i_1})\ldots 
f_p(\sigma_{i_p}) \Bigr\rangle^n.
\end{eqnarray*}
Using replicas $\vsi^1,\ldots, \vsi^n$ we have
$$\Bigl\langle f_1(\sigma_{i_1})\ldots f_p(\sigma_{i_p}) \Bigr\rangle^n =
\Bigl\langle \prod_{l\leq n} f_1(\sigma_{i_1}^l)\ldots f_p(\sigma_{i_p}^l) \Bigr\rangle
$$
and, thus,
$$
\frac{1}{N^p}\sum_{i_1,\ldots,i_p=1}^{N} \Bigl\langle f_1(\sigma_{i_1})\ldots 
f_p(\sigma_{i_p}) \Bigr\rangle^n =
\Bigl\langle \prod_{j\leq p} A_{j,n} \Bigr\rangle,
$$
where
$$
A_{j,n}=A_{j,n}(\vsi^1,\ldots, \vsi^n) = \frac{1}{N}\sum_{i\leq N} \prod_{l\leq n}f_j(\sigma_i^l).
$$
Denote by $\e_0$ the expectation in $f_1,\ldots,f_p$ 
and $x_1,\ldots,x_p$ only (the r.v. $x_1,\ldots,x_p$ 
are not present here and will appear in the terms below).
Since
$f_1,\ldots,f_p$ are i.i.d. and independent of the randomness in 
$\la\cdot\ra,$
$\e_0 \la\prod_{j\leq p} A_{j,n}\ra =  \la \e_0 \prod_{j\leq p} A_{j,n}\ra=\la B_n^p\ra$ where $B_n = \e_0 A_{j,n}.$
Therefore, since we also assumed that $b$ is independent of $f_1,\ldots,f_p,$
\begin{equation}
\e_0 \frac{1}{N^p}
\sum_{i_1,\ldots,i_p=1}^N
\log \Bigl\langle \exp \theta(\sigma_{i_{1}},\ldots,\sigma_{i_{p}})\Bigr\rangle
= \e_0 \log a 
- \sum_{n=1}^{\infty} \frac{(-b)^n}{n} \Bigl\langle B_n^p \Bigr\rangle.
\label{neg1}
\end{equation}
A similar analysis applies to the second term in (\ref{neg}).
First of all, (\ref{E1}) implies that
\begin{eqnarray*}
\exp U(\sigma_i) = 
\la \exp \theta \ra_{\vec{x}}^{-} \Bigr|_{\eps_p=\sigma_i}
=
a\Bigl(1+b f_p(\sigma_i)\prod_{l=1}^{p-1} 
\frac{\Av f_l(\eps)\exp \eps x_l}
{\ch (x_l)}\Bigr),
\end{eqnarray*}
and, therefore,
\begin{eqnarray*}
&&
\log \Bigl\langle \exp U(\sigma_i)\Bigr\rangle = \log a 
-\sum_{n=1}^{\infty}\frac{(-b)^n}{n}
\Bigl(
\Bigl\langle f_p(\sigma_i)\Bigr\rangle \prod_{l=1}^{p-1} 
\frac{\Av f_l(\eps)\exp \eps x_l} {\ch (x_l)}
\Bigr)^n
\\
&&
=
\log a 
-\sum_{n=1}^{\infty}\frac{(-b)^n}{n}
\Bigl\langle f_p(\sigma_i^1)\ldots f_p(\sigma_i^n)\Bigr\rangle 
\prod_{l=1}^{p-1} \Bigl(
\frac{\Av f_l(\eps)\exp \eps x_l} {\ch (x_l)}
\Bigr)^n,
\end{eqnarray*}
where in the last equality we used replicas. Now,
\begin{eqnarray*}
\frac{1}{N} \sum_{i=1}^{N}\log \Bigl\langle \exp U(\sigma_i)\Bigr\rangle
=
\log a 
-\sum_{n=1}^{\infty}\frac{(-b)^n}{n}
\Bigl\langle A_{p,n}\Bigr\rangle 
\prod_{l=1}^{p-1} \Bigl(
\frac{\Av f_l(\eps)\exp \eps x_l} {\ch (x_l)}
\Bigr)^n
\end{eqnarray*}
and taking the expectation w.r.t. $f_1,\ldots,f_p$ 
and $x_1,\ldots,x_p$ we get
\begin{equation}
\e_0 \frac{1}{N} \sum_{i=1}^{N}\log \Bigl\langle \exp U(\sigma_i)\Bigr\rangle
=
\e_0 \log a 
-\sum_{n=1}^{\infty}\frac{(-b)^n}{n}
\Bigl\langle B_n\Bigr\rangle (C_n)^{p-1}
\label{neg2}
\end{equation}
where
$$
C_n = \e_0
\Bigl(
\frac{\Av f_l(\eps)\exp \eps x_l} {\ch (x_l)}
\Bigr)^n.
$$
Finally, in absolutely similar manner
\begin{equation}
\e_0 \log \la \E \ra_{\vec{x}} =
\e_0 \log a 
-\sum_{n=1}^{\infty}\frac{(-b)^n}{n}
(C_n)^p.
\label{neg3}
\end{equation}
Combining (\ref{neg1}), (\ref{neg2}) and (\ref{neg3}) we see that  (\ref{neg}) can be written as
\begin{equation}
-\sum_{n=1}^{\infty}\frac{\e (-b)^n}{n}
\e \Bigl\langle
B_n^p -p B_n C_n^{p-1} + (p-1) (C_n)^p \Bigr\rangle \leq 0
\label{negdone}
\end{equation}
which holds true using condition (\ref{condition2}) and the fact that
$x^p - pxy^{p-1} +(p-1)y^p\geq 0$ for all $x,y\in\Reals.$
This finishes the proof of Theorem \ref{RS}.
\qed

\section{A general weighting scheme.}

The use of weighting scheme as considered in this section is directly motivated by the paper \cite{ASS}. It is a very useful device, see e.g. \cite{T2}. 

We consider a countable index set $\G$,  an arbitrary sequence of r.v. 
$(x^{\gamma})_{\gamma\in\G}$ and let $(x_l^{i,j,\gamma})_{\gamma\in\G}$
for $i,j,l\geq 1$ be its independent copies of this sequence. 
We define
$$
U_{i,j}^{\gamma}(\eps)=
U(\theta_{i,j}, x_1^{i,j,\gamma},\ldots,x_{p-1}^{i,j,\gamma},\eps)
$$
where $\theta_{i,j}$ are independent copies of $\theta$
and consider the Hamiltonian
\begin{equation} 
-H_{N,t}^{\gamma}(\vsi)=\sum_{k\leq M_t}\theta_k(\sigma_{i_{1,k}},
\ldots,\sigma_{i_{p,k}}) + 
\sum_{i\leq N}\sum_{j\leq k_{i,t}} U_{i,j}^{\gamma}(\sigma_i) 
+H'_N(\vsi),
\label{inter-weights}
\end{equation}
where $M_t$ and $k_{i,t}, i\leq N$ are defined as in (\ref{inter-rs}).

Consider an arbitrary random sequence $(v_\gamma)_{\gamma\in\G}$
independent of all r.v. in (\ref{inter-weights})
and such that $\sum_{\gamma\in\G}v_{\gamma}=1$ and define
the Gibbs measure on $\Sigma_N\times \G$ by
\begin{equation}
G(\{\vsi,\gamma\})=
v_{\gamma} \exp( - H_{N,t}^{\gamma}(\vsi))/Z_N
\nonumber
\end{equation}
where the partition function $Z_N$ is given by 
$
Z_N = \sum_{\gamma,\vsi} v_{\gamma} \exp(-H_{N,t}^{\gamma}(\vsi)).
$
For a function
$f(\vsi,\gamma)$ on $\Sigma_N\times\Gamma,$ 
$\la\cdot\ra$ will now denote the average
with respect to the Gibbs measure
$G$
\begin{equation}
\la f \ra = 
\frac{1}{Z_N}
\sum_{\gamma\in\G,\vsi\in\Sigma_N} 
f(\vsi,\gamma) v_{\gamma} \exp(-H_{N,t}^{\gamma}(\vsi)).
\label{WGibbs}
\end{equation}
Let
$$
\varphi(t)=\frac{1}{N}\e \log 
\sum_{\gamma\in\G}\sum_{\vsi\in\Sigma_N} 
v_{\gamma} \exp(-H_{N,t}^{\gamma}(\vsi)).
$$
Clearly, $F_N = \varphi(1).$ The following Theorem holds.

\begin{theorem}\label{W}.
If conditions (\ref{condition1}),
(\ref{condition2}) and (\ref{condition3}) hold then
\begin{equation}
F_N \leq \varphi(0) - \alpha(p-1) 
\int_{0}^{1}
\e\log \Bigl\langle\la \E 
\ra_{\vec{x}^{\gamma}}\Bigr\rangle dt.
\label{Wbound}
\end{equation}
 where $\vec{x}^{\gamma}=(x_1^{\gamma},\ldots,x_p^{\gamma})$ 
and where $(x_l^{\gamma})_{\gamma\in\G}$
are independent copies of $(x^{\gamma})_{\gamma\in\G}$
for $l\leq p.$
\end{theorem}
Of course, the integrand $\e\log\la\la\E\ra_{\vec{x}^{\gamma}}\ra$
in the last term of (\ref{Wbound}) depends on $t$ through
$\la\cdot\ra$ since $\la\E\ra_{\vec{x}^{\gamma}}$ is a function
of $\gamma.$

{\bf Proof.}
The proof follows that of Theorem \ref{RS} with almost no changes.
(\ref{deriv}) now becomes
\begin{equation}
\varphi'(t)=\alpha\Bigl(
\frac{1}{N^p}
\sum_{i_1,\ldots,i_p=1}^N
\e \log 
\Bigl\langle \exp \theta(\sigma_{i_{1}},\ldots,\sigma_{i_{p}})\Bigr\rangle
-p\frac{1}{N}\sum_{i=1}^N
\e \log \Bigl\langle \exp U^{\gamma}(\sigma_i) \Bigr\rangle
\Bigr)
\label{Wderiv}
\end{equation}
where $\la\cdot\ra$ is given by (\ref{WGibbs}).
Similarly to (\ref{neg}) we will now show that
\begin{equation}
\frac{1}{N^p}
\sum_{i_1,\ldots,i_p=1}^N
\e \log 
\Bigl\langle \exp \theta(\sigma_{i_{1}},\ldots,\sigma_{i_{p}})\Bigr\rangle
-p\frac{1}{N}\sum_{i=1}^N
\e \log \Bigl\langle \exp U^{\gamma}(\sigma_i) \Bigr\rangle
+(p-1)\e \log \Bigl\langle \la \E \ra_{\vec{x}^{\gamma}} \Bigr\rangle \leq 0.
\label{Wneg}
\end{equation}
This clearly implies the statement of Theorem \ref{W}
since if we denote $c(t)=\alpha(p-1)\e\log\la\la\E\ra_{\vec{x}^{\gamma}}\ra$
then equation (\ref{Wneg}) yields
$\varphi'(t)+c(t)\leq 0$ and therefore
$$
\varphi(1)\leq \varphi(0) - \int_{0}^{1} c(t) dt
$$
which is precisely the statement of the Theorem.

In the proof of Theorem \ref{RS} we showed that (\ref{neg})
is equivalent to (\ref{negdone}).
Following the same arguments one can show that
(\ref{Wneg}) can be written as
$$
-\sum_{n=1}^{\infty}\frac{\e (-b)^n}{n}
\e \Bigl\langle
B_n^p -p B_n C_n^{p-1}(\gamma_1,\ldots,\gamma_n) + 
(p-1) C_n^p(\gamma_1,\ldots,\gamma_n) \Bigr\rangle \leq 0
$$
where now
$$
C_n(\gamma_1,\ldots,\gamma_n)=
\e_0
\prod_{j=1}^{n}
\frac{\Av f_1(\eps)\exp \eps x_1^{\gamma_j}} {\ch (x_1^{\gamma_j})}.
$$
The one difference with the case of Theorem \ref{RS} is  the fact that using replicas to represent $\la\cdot\ra^n$
now involves both $\vsi$ and $\gamma.$
For instance, in the calculations leading to (\ref{neg2}) 
there will appear a term
$$
\Bigl\langle f_p(\sigma_i) \prod_{l=1}^{p-1} 
\frac{\Av f_l(\eps)\exp \eps x_l^{\gamma}} {\ch (x_l^{\gamma})} \Bigr\rangle^n
$$
where $(x_l^{\gamma})_{\gamma\in\G}, l\geq 1$ are independent copies of
$(x^{\gamma})_{\gamma\in\G}$ independent of the randomness
in $\la\cdot\ra.$
Using replicas $(\vsi^1,\gamma_1),\ldots,(\vsi^n,\gamma_n)$, 
this term can be written as
$$
\Bigl\langle
f_p(\sigma_i^1)\ldots f_p(\sigma_i^n) 
\prod_{j=1}^{n}\prod_{l=1}^{p-1} 
\frac{\Av f_l(\eps)\exp \eps x_l^{\gamma_j}} {\ch (x_l^{\gamma_j})}
\Bigr\rangle.
$$
Averaging for $i$ and taking expectation $\e_0$ w.r.t.
$f_1,\ldots,f_p$ and $(x_l^{\gamma})$
yields
$$
\e_0 \Bigl\langle
A_{p,n} 
\prod_{l=1}^{p-1}
\prod_{j=1}^{n}
\frac{\Av f_l(\eps)\exp \eps x_l^{\gamma_j}} {\ch (x_l^{\gamma_j})}
\Bigr\rangle
=
\Bigl\langle
B_n C_n^{p-1}(\gamma_1,\ldots,\gamma_j)
\Bigr\rangle. 
$$
Similarly, the calculations leading to (\ref{neg3})
will produce $\la C_n^p(\gamma_1,\ldots,\gamma_p) \ra.$
The rest of the argument is the same.

\qed

\section {The $\vec{1}$-step of replica-symmetry breaking bound.}

In the context and with the notations of the previous section
we will now make specific choices of
the random sequences $(v_{\gamma})$
and $(x_l^{i,j,\gamma}).$  From now on we will also assume that
$H'_N(\vsi)$ is given by (\ref{H2}).

We denote by  $\M_1$ the  set of probability measures on $\Reals,$
and $\M_2$ the  set of probability measures on $\M_1.$
Consider  $\zeta\in\M_2$, our basic parameter on which will depend the bound we are going to obtain.
We consider a sequence $(\eta_l, x_l)_{l \geq 1}$ with the following properties. The sequence $(\eta_l)$ is an i.i.d sequence of $\M_1$ distributed according to $\zeta$. Conditionally on this sequence, the sequence $(x_l)$ is independent and $x_l$ is distributed like $\eta_l$. We consider i.i.d. copies $(\eta_{l}^j, x_{l}^j)$ of the sequence $(\eta_l, x_l)$.

\begin{theorem}\label{1step}
Suppose that (\ref{condition1}), (\ref{condition2}) and
(\ref{condition3}) hold  and $H'_N(\vsi)$ is given by (\ref{H2}). 
Let $U_j(\eps)=U(\theta_j, x_1^j,\ldots,x_{p-1}^j, \eps),$
where $\theta_j$ are independent copies of $\theta,$ 
$\vec{x}=(x_1,\ldots,x_{p})$ and let $k$ be a Poisson r.v.
with mean $\e k = \alpha p.$ Then, for $m\in (0,1)$ we have
\begin{equation}
F_N\leq \Phi_1(\zeta,m)=
\log 2 + \frac{1}{m} \e\log\e'\Bigl(
\Av\exp\bigl(\sum_{j\leq k} U_j(\eps) + h\eps\bigr)\Bigr)^m
-\alpha(p-1)\frac{1}{m}\e\log\e'(\la\E\ra_{\vec{x}})^m, 
\label{1RSB}
\end{equation}
where 
$\e'$ is the expectation w.r.t. $(x_l)$ and $(x_l^j)$ 
for fixed $(\eta_l)$ and $(\eta_l^j)$ and 
$\e$ denotes the expectation w.r.t. 
$(\eta_l),$ $(\eta_l^j),$ $(\theta_j),$ $k$ and $h.$  
\end{theorem}

Of course, Theorem \ref{1step} implies that
$$
F_N\leq \inf_{\zeta,m}\Phi_1(\zeta,m).
$$
It should also be noted that Theorem~\ref{1step} is a generalization of Theorem~\ref{RS}, as is seen by taking $m \rightarrow 0$ and $\zeta$ concentrated at one point of $\M_1$.

The terms $\e'(\cdot)^m$ will magically appear with the proper choice of the sequence $v_\gamma$, that we explain first. Let $\G$ be a set of natural numbers $\Natural.$
We consider a non-increasing enumeration  $(u_{\gamma})_{\gamma\geq 1}$ of the points 
generated by Poisson process on $\Reals^{+}$ of intensity
measure $x^{-m-1}.$ To avoid repetition, we will say that such a sequence has distribution $\Xi_m$.
We define a sequence $(v_{\gamma})_{\gamma\geq 1}$
by
$$
v_{\gamma}=\frac{u_{\gamma}}{\sum_{\gamma'\in\G} u_{\gamma'}}.
$$
%It is well defined (see \cite{SG}, Section 1.2). 
With the notation of \cite{SG}, this sequence has the Poisson-Dirichlet distribution $\Lambda_m$.
This key property is as follows 
(see e.g. Proposition 6.5.15 in \cite{SG}).
\begin{proposition}\label{PD}
Consider a r.v. $\xi\geq 0,$  $\e\xi^2<\infty$ and independent copies
$(\xi_\gamma)_{\gamma\geq 1}.$ Then the sequences
$(u_{\gamma} \xi_{\gamma})_{\gamma\geq 1}$ and 
$\bigl(u_\gamma (\e \xi_1^m)^{1/m} \bigr)_{\gamma\geq 1}$ 
have the same distribution and, therefore,
\begin{equation}
\e\log \sum_{\gamma\geq 1} v_{\gamma} \xi_{\gamma}=\e \log \sum_{\gamma \geq 1}u_\gamma \xi_\gamma -\e \log \sum_{\gamma \geq 1} u_\gamma
=\frac{1}{m}\log \e \xi^m.
\label{PDprop}
\end{equation}
\end{proposition}

{\bf Proof of Theorem \ref{1step}.}
We consider an element $\eta$ of $\M_1$ that is distributed according to $\zeta$ and a sequence $(x^\gamma)_{\gamma \geq 1}$ that, given $\eta$, is i.i.d distributed according to $\eta$. For $i,j,l\geq 1$ we consider independent copies $(\eta_l^{i,j})$ and $(x_l^{i,j,\gamma})$ of these variables. We also consider other independent copies $(\eta_l)$ and $(x_l^\gamma)$ of these variables. We denote by $\F$ the $\sigma$-algebra generated by the variables $\eta_l, \eta_l^{i,j}, h_i, k_{i,j}$ and $\theta_{i,j}$.
%Let us also consider sequences $(x_l)$ and $(x_l^j)$ that are independent copies of $(x_l^{1,1,1})$ and $(x_l^{j,1,1})$ correspondingly.
%{\bf Remark.} The bound (\ref{1RSB}) is commonly written in terms of functions
%$B$ and $u$ defined in (\ref{Bu}). Namely, one can write
%$$
%\Av\exp\sum_{j\leq k}U_j(\eps)=
%\Av\prod_{j\leq k}\la \E \ra_{\vec{x^{j}}}^{-}\bigr|_{\theta=\theta_j}
%=
%\Av\prod_{j\leq k}B_j \exp( \eps u_j)
%=\ch\bigl(\sum_{j\leq k} u_j\bigr)\prod_{j\leq k}B_j.
%$$

Let us first consider the integrand
$\e\log \bigl\langle\la \E \ra_{\vec{x}^{\gamma}}\bigr\rangle$
in the last term of (\ref{Wbound}).
Let us denote 
$$
Z_{t}(\gamma)=\sum_{\vsi} \exp(-H_{N,t}^{\gamma}(\vsi))
$$
and $e(\gamma) = \la \E \ra_{\vec{x}^{\gamma}}.$
Note that
$e(\gamma)$ depends on $x_l^{\gamma}, l\geq 1$ and $Z_t(\gamma)$ depends
on $x_l^{i,j,\gamma}, i,j,l\geq 1$ through the Hamiltonian $H_{N,t}^{\gamma}.$
By construction, given  $\F$, the
 sequences $(Z_{t}(\gamma))_{\gamma \geq 1}$ and $(e(\gamma))_{\gamma\geq 1}$ 
are i.i.d. and independent of each other.
If we denote by $\e'$ conditional  expectation given $\F$ then,
using Proposition \ref{PD}, we get
\begin{eqnarray*}
&&
\e'\log \Bigl\langle\la \E \ra_{\vec{x}^{\gamma}}\Bigr\rangle
=
\e' \log 
\frac{\sum v_{\gamma} e(\gamma) Z_{t}(\gamma)}
{\sum v_{\gamma}Z_{t}(\gamma)}
=
\e' \log 
\sum v_{\gamma} e(\gamma) Z_{t}(\gamma) -
\e' \log 
\sum v_{\gamma}Z_{t}(\gamma)
\\
&&
=
\frac{1}{m}\log \e' \bigl(e Z_{t}\bigr)^m - 
\frac{1}{m} \log \e' Z_{t}^m
=
\frac{1}{m}\log \bigl(\e' e^m \e' Z_{t}^m \bigr) - 
\frac{1}{m} \log \e' Z_{t}^m =
\frac{1}{m}\log \e' e^m,
\end{eqnarray*}
which is independent of $t$ and, therefore,
this yields the last term of (\ref{1RSB}).

Next let us consider the first term in (\ref{Wbound}), $\varphi(0).$
First of all,
$$
\varphi(0)=\frac{1}{N}\e\log\sum_{\gamma\geq 1}v_{\gamma} Z_0(\gamma),
$$
where
\begin{eqnarray*}
Z_0(\gamma)=\sum_{\vsi} \exp(-H_{N,0}^{\gamma}(\vsi))
=2^N\prod_{i=1}^{N}\Av\exp\Bigl(
\sum_{j\leq k_{i,0}} U_{i,j}^{\gamma}(\eps) + h_i\eps\Bigr)
=2^N\prod_{i=1}^{N}F_i(\gamma),
\end{eqnarray*}
and where we introduced the notation 
$F_i(\gamma)= \Av \exp\Bigl(\sum_{j\leq k_{i,0}} 
U_{i,j}^{\gamma}(\eps) +h_i\eps\Bigr).$
By construction, given $\F$, the sequences $F_i(\gamma), \gamma\geq 1$
are i.i.d. and independent for different indices $i.$
Therefore, application of Proposition \ref{PD} gives
\begin{eqnarray*}
&&
\frac{1}{N}\e'\log 2^N\sum_{\gamma\geq 1}v_{\gamma}\prod_{i=1}^{N}F_i(\gamma)
=\log 2 + \frac{1}{N}\frac{1}{m}\log\e'\Bigl(\prod_{i=1}^{N}F_i(1)\Bigr)^m
\\
&&
=\log 2 + \frac{1}{N}\sum_{i=1}^{N}\frac{1}{m}\log\e'\bigl(F_i(1)\bigr)^m,
\end{eqnarray*}
and taking expectation w.r.t. all the other r.v. implies
$$
\varphi(0)=\log 2 + \frac{1}{m}\e\log\e'\bigl(F_1(1)\bigr)^m,
$$
which is precisely the first two terms in (\ref{1RSB}).

\qed

\section{The $\vec{r}$-step of replica-symmetry breaking bound.}

We will first explain the choice of weights $(v_{\gamma})$
and r.v. $(x^{\gamma})$ in Theorem \ref{W} that will yield
the $r$-step of replica symmetry breaking bound of Theorem \ref{rstep}
below.

We take $\G=\Natural^r$ and
define a sequence $(v_\gamma)_{\gamma\in\G}$ 
using Derrida-Ruelle cascades (see \cite{Ruelle}).
Consider arbitrary parameters $0<m_1<\ldots<m_r<1.$
Let us consider sequences $(u_{\gamma_1})_{\gamma_1\geq 1},
\ldots, (u_{\gamma_r})_{\gamma_r\geq 1}$ with the distributions
$\Xi_{m_1},\ldots,\Xi_{m_r}$ correspondingly.
For $2\leq l\leq r$ let us consider a sequence
$(u_{\gamma_1,\ldots,\gamma_l})_{\gamma_1,\ldots,\gamma_l\geq 1}$ 
that for any fixed $(\gamma_1,\ldots,\gamma_{l-1})$ is
an independent copy of the sequence $(u_{\gamma_l})_{\gamma_l\geq 1}.$ 
We define $\bar{u}_{\gamma_1,\ldots,\gamma_r}=\prod_{l=1}^{r}
u_{\gamma_1,\ldots,\gamma_l}$ and
\begin{equation}
v_{\gamma_1,\ldots,\gamma_r}=\frac{\bar{u}_{\gamma_1,\ldots,\gamma_r}}
{\sum_{\gamma_1',\ldots,\gamma_r'}
\bar{u}_{\gamma_1',\ldots,\gamma_r'}}.
\label{rweights}
\end{equation}
Next we define the set of r.v. $x_{\omega}(\gamma_1,\ldots,\gamma_r)$
for $\omega\in\Omega$ and $\gamma_1,\ldots,\gamma_r\geq 1.$

Let $\M_1$ be a set of probability measures on $\Reals,$ and
by induction for $l\leq r$ we define $\M_{l+1}$ as a set
of probability measures on $\M_{l}.$
Let us fix $\zeta \in \M_{r+1}$ (our basic parameter) and define a random sequence 
$(\eta, \eta(\gamma_1), \ldots, \eta(\gamma_1, \ldots, \gamma_{r-1}), x(\gamma_1,\ldots, \gamma_r))$ as follows.  
The element $\eta$ of $\M_r$ is distributed according to $\zeta$. Given $\eta$, 
the sequence $(\eta(\gamma_1))_{\gamma_1 \geq 1}$ of elements of $\M_{r-1}$ is i.i.d distributed like $\eta$.  
For $1\leq l\leq r-1$,  given  all the elements 
$\eta(a_1, \ldots, a_{s})$ for all values of the integers $a_1, \ldots, a_s$ and all $s \leq l-1$,
 the sequence $(\eta(\gamma_1, \ldots, \gamma_l))_{\gamma_l \geq 1}$ of elements of $\M_{r-l}$  
is i.i.d  distributed like $\eta(\gamma_1, \ldots, \gamma_{l-1})$, and these sequences are independent 
of each other for different values of $ (\gamma_1, \ldots, \gamma_{l-1})$.   
Finally, given  all the elements $\eta(a_1, \ldots, a_{s})$ for all values of the integers 
$a_1, \ldots, a_s$ and all $s \leq r-1$ the sequences  
$x(\gamma_1,\ldots,\gamma_r), \gamma_r\geq 1$ 
is  an i.i.d. sequence on $\Reals$ with the distribution
$\eta(\gamma_1,\ldots,\gamma_{r-1})$ and these sequences
are independent for different values of  $(\gamma_1,\ldots,\gamma_{r-1}).$
The process of generating $x$'s can be represented schematically as
\begin{equation}
\zeta\to\eta \to\eta(\gamma_1)\to\ldots
\to\eta (\gamma_1,\ldots,\gamma_{r-1})\to
x(\gamma_1,\ldots,\gamma_r).
\label{xs}
\end{equation}

For simplicity of notations instead of writing various combination
of indices $i,j,l$ let us first consider an arbitrary countable index
set $\Omega.$ For $\omega \in \Omega$, we consider independent copies 
$(\eta_\omega, \eta_\omega(\gamma_1), \ldots, \eta_\omega(\gamma_1, \ldots, 
\gamma_{r-1}), x_\omega(\gamma_1,\ldots, \gamma_r))$ of $(\eta, 
\eta(\gamma_1), \ldots,  $ $\eta(\gamma_1, \ldots, \gamma_{r-1}), 
x(\gamma_1,\ldots, \gamma_r))$.

For $ 0\leq j\leq r-1,$
let us denote by $\F_j$ the $\sigma$-algebra generated by 
$\eta_{\omega}(\gamma_1,\ldots,\gamma_l)$
for $\omega\in\Omega,$ $l\leq j,$ $\gamma_1,\ldots,\gamma_l\geq 1$, and by the r.v $h_i$, 
$\theta_{i,j}$ and $k_{i,j}$.
Let us denote by $\e_{j}$ the expectation given $\F_j$ or, in other words,
w.r.t. $\eta_{\omega}(\gamma_1,\ldots,\gamma_l)$
for $\omega\in\Omega,$ $l>j,$ $\gamma_1,\ldots,\gamma_l\geq 1$
and
$x_{\omega}(\gamma_1,\ldots,\gamma_r)$
for $\omega\in\Omega,$ $\gamma_1,\ldots,\gamma_r\geq 1.$ In particular $\F_0$ is 
generated by the variables $\eta_\omega$, $h_i$, $\theta_{i,j}$ and $k_{i,j}$.

For a random variable $U \geq 0$ we define $T_rU = U$ and by induction, for $0 \leq l<r$ 
we define the r.v.
$U_l $  by
\begin{equation}
T_l U = \Bigl(
\e_{l} (T_{l+1} U)^{m_{l+1}}
\Bigr)^{1/m_{l+1}}.
\label{T}
\end{equation}
Let us consider a function $V:\Reals^{\Omega}\to\Reals$, $V\geq 0$, 
and the r.v.
\begin{equation}
V(\gamma_1,\ldots,\gamma_r)=
V\Bigl(\bigl(x_{\omega}(\gamma_1,\ldots,\gamma_r)\bigr)_{\omega\in\Omega}
\Bigr).
\label{Vee}
\end{equation}
 The distribution of the r.v. $V(\gamma_1, \ldots, \gamma_r)$ is independent 
of the value $(\gamma_1, \ldots, \gamma_r)$, and the r.v. $T_l(V(\gamma_1, 
\ldots, \gamma_r))$ depends only on $\gamma_1, \ldots, \gamma_l$. 

The following key property is based on iterative application 
of Proposition \ref{PD}. It should be obvious to a reader familiar
with Derrida-Ruelle cascades (\cite{Ruelle}). In fact the essential ideas of the present scheme of proof are  apparently  known to the authors of \cite{ASS}, but for lack of references, it seems appropriate to give complete details. 

\begin{proposition}\label{PD2}
If $V$ is defined by (\ref{Vee}) and
$\e V^2<\infty$ then
\begin{equation}
\e\log\sum_{\gamma_1,\ldots,\gamma_r\geq 1}
v_{\gamma_1,\ldots,\gamma_r}V(\gamma_1,\ldots,\gamma_r)
=\e\log T_0 V.
\label{PDprop2}
\end{equation}
\end{proposition}
{\bf Proof.}
Let $\e'$ denote the expectation w.r.t. $(v_{\gamma_1,\ldots,\gamma_r})$
and $(x_{\omega}(\gamma_1,\ldots,\gamma_r))$ given $\F_0,$
i.e. for a fixed sequence $(\eta_{\omega}).$
Let us first consider
$$
\e'\log\sum_{\gamma_1,\ldots,\gamma_r\geq 1}
\bar{u}_{\gamma_1,\ldots,\gamma_r}V(\gamma_1,\ldots,\gamma_r).
$$
By the definition of $\bar{u}_{\gamma_1,\ldots,\gamma_{r}}$ we can write
\begin{eqnarray*}
\sum_{\gamma_1,\ldots,\gamma_r\geq 1}
\bar{u}_{\gamma_1,\ldots,\gamma_r}V(\gamma_1,\ldots,\gamma_r)=
\sum_{\gamma_1,\ldots,\gamma_{r-1}\geq 1}
\prod_{l\leq r-1} u_{\gamma_1,\ldots,\gamma_l}
\Bigl(
\sum_{\gamma_r\geq 1}u_{\gamma_1,\ldots,\gamma_{r}}V(\gamma_1,\ldots,\gamma_r)
\Bigr).
\end{eqnarray*}
For a fixed $(\gamma_1,\ldots,\gamma_{r-1})$, and 
given $\F_{r-1}$,
the sequence 
$V(\gamma_1,\ldots,\gamma_r),\gamma_r\geq 1$ is i.i.d.
while the sequence  $(u_{\gamma_1,\ldots,\gamma_r})_{\gamma_r\geq 1}$ has distribution
$\Xi_{m_r}.$ Therefore, writing for simplicity $T_{r-1}V$ rather than $T_{r-1}V(\gamma_1, \cdots, \gamma_r)$, Proposition \ref{PD} implies that
\begin{equation}
\sum_{\gamma_r\geq 1}\bar{u}_{\gamma_1,\ldots,\gamma_{r}}V(\gamma_1,\ldots,\gamma_r)
\sim
(T_{r-1}V) \sum_{\gamma_r\geq 1}u_{\gamma_1,\ldots,\gamma_{r}}
= (T_{r-1}V)S_{r-1}(\gamma_1,\ldots,\gamma_{r-1}),
\label{pd1}
\end{equation}
where $\sim $ means equality in distribution and where
we introduced the notation
$$
S_{r-1}(\gamma_1,\ldots,\gamma_{r-1})=\sum_{\gamma_r\geq 1}u_{\gamma_1,\ldots,\gamma_{r-1},\gamma_{r}}.
$$
Of course, $T_{r-1}V$  depends on $(\gamma_1,\ldots,\gamma_{r-1})$, although this is not explicit in the notation.

Moreover, given $\F_{r-1},$ both sides of (\ref{pd1}) 
are by construction independent for different indices $(\gamma_1,\ldots,\gamma_{r-1})$
and, thus,
\begin{eqnarray}
&&
\sum_{\gamma_1,\ldots,\gamma_r\geq 1}
\bar{u}_{\gamma_1,\ldots,\gamma_r}V(\gamma_1,\ldots,\gamma_r)
\sim
\sum_{\gamma_1,\ldots,\gamma_{r-1}\geq 1}
\prod_{l\leq r-1} u_{\gamma_1,\ldots,\gamma_l}
(T_{r-1}V)S_{r-1}(\gamma_1,\ldots,\gamma_{r-1})
\nonumber
\\
&&
=
\sum_{\gamma_1,\ldots,\gamma_{r-2}\geq 1}
\prod_{l\leq r-2} u_{\gamma_1,\ldots,\gamma_l}
\Bigl(
\sum_{\gamma_{r-1}\geq 1} u_{\gamma_1,\ldots,\gamma_{r-1}}
(T_{r-1}V)S_{r-1}(\gamma_1,\ldots,\gamma_{r-1})\Bigr).
\label{s1}
\end{eqnarray}
For a fixed $(\gamma_1,\ldots,\gamma_{r-2})$, and 
given $\F_{r-2}$,
%fixed $(\eta_{\omega}(\gamma_1,\ldots,\gamma_{r-1}))_{\omega\in\Omega}$
the sequences 
$$
(T_{r-1}V)(\gamma_1,\ldots,\gamma_{r-1}) \mbox{ and }
S_{r-1}(\gamma_1,\ldots,\gamma_{r-1}) \mbox{ for } 
\gamma_{r-1}\geq 1
$$ 
are i.i.d.
and by construction the sequence $(u_{\gamma_1,\ldots,\gamma_{r-1}})_{\gamma_{r-1}\geq 1}$ has distribution
$\Xi_{m_{r-1}}.$ Therefore, Proposition \ref{PD} now implies that
\begin{eqnarray}
&&
\sum_{\gamma_{r-1}\geq 1}u_{\gamma_1,\ldots,\gamma_{r-1}}(T_{r-1}V)(\gamma_1,\ldots,\gamma_{r-1})
S_{r-1}(\gamma_1,\ldots,\gamma_{r-1})
\sim
(T_{r-2}V) C_{m_{r},m_{r-1}}
\sum_{\gamma_{r-1}\geq 1}u_{\gamma_1,\ldots,\gamma_{r-1}}
\nonumber
\\
&&
= (T_{r-2}V) C_{m_{r},m_{r_1}}
S_{r-2}(\gamma_1,\ldots,\gamma_{r-2}),
\label{pd2}
\end{eqnarray}
where
$S_{r-2}(\gamma_1,\ldots,\gamma_{r-2})=
\sum_{\gamma_{r-1}\geq 1}u_{\gamma_1,\ldots,\gamma_{r-1}}$ and
$$
C_{m_r,m_{r-1}}=\Bigl(\e
\bigl(S_{r-1}(\gamma_1,\ldots,\gamma_{r-1})\bigr)^{m_{r-1}}
\Bigr)^{1/m_{r-1}}
=
\Bigl(\e
\bigl(\sum_{\gamma_r\geq 1} u_{\gamma_r}\bigr)^{m_{r-1}}
\Bigr)^{1/m_{r-1}},
$$
where $(u_{\gamma_r})_{\gamma_r\geq 1}$ has distribution $\Xi_{m_{r}}.$
One can easily check that $C_{m_r,m_{r-1}}<\infty$
due to the fact that $m_{r-1}< m_{r}.$

Given $\F_{r-2},$ both sides of (\ref{pd2}) 
are independent for different $(\gamma_1,\ldots,\gamma_{r-2})$
and, therefore, (\ref{s1}) implies
\begin{eqnarray*}
&&
\sum_{\gamma_1,\ldots,\gamma_r\geq 1}
\bar{u}_{\gamma_1,\ldots,\gamma_r}V(\gamma_1,\ldots,\gamma_r)
\sim
\sum_{\gamma_1,\ldots,\gamma_{r-2}\geq 1}
\prod_{l\leq r-2} u_{\gamma_1,\ldots,\gamma_l}
(T_{r-2}V) C_{m_{r},m_{r_1}}
S_{r-2}(\gamma_1,\ldots,\gamma_{r-2})
\\
&&
=
\sum_{\gamma_1,\ldots,\gamma_{r-3}\geq 1}
\prod_{l\leq r-3} u_{\gamma_1,\ldots,\gamma_l}
\Bigl(
\sum_{\gamma_{r-2}\geq 1} u_{\gamma_1,\ldots,\gamma_{r-2}}
(T_{r-2}V) C_{m_{r},m_{r-1}}
S_{r-2}(\gamma_1,\ldots,\gamma_{r-2})\Bigr).
\end{eqnarray*}
Let us define
$$
C_{m_j,m_{j-1}}
=
\Bigl(\e
\bigl(\sum_{\gamma_j\geq 1} u_{\gamma_j}\bigr)^{m_{j-1}}
\Bigr)^{1/m_{j-1}},
$$
where the sequence $(u_{\gamma_j})_{\gamma_j\geq 1}$ has distribution $\Xi_{m_{j}}$
and let 
$$
C_{j-1}= C_{m_r,m_{r-1}}\cdots C_{m_j,m_{j-1}}.
$$
Denote
$$
S_{j}(\gamma_1,\ldots,\gamma_{j})=\sum_{\gamma_{j+1}\geq 1}u_{\gamma_1,\ldots,\gamma_{j+1}}.
$$
Repeating the same argument as above one can show 
by decreasing induction over $j$ that
\begin{equation}
\sum_{\gamma_1,\ldots,\gamma_r\geq 1}
\bar{u}_{\gamma_1,\ldots,\gamma_r}V(\gamma_1,\ldots,\gamma_r)
\sim
\sum_{\gamma_1,\ldots,\gamma_{j}\geq 1}
\prod_{l\leq j} u_{\gamma_1,\ldots,\gamma_j}
\Bigl(
(T_{j}V) C_{j+1}
S_{j}(\gamma_1,\ldots,\gamma_{j})\Bigr).
\label{sj}
\end{equation}
In particular, for $j=0,$ (\ref{sj}) reads
\begin{eqnarray*}
\sum_{\gamma_1,\ldots,\gamma_r\geq 1}
\bar{u}_{\gamma_1,\ldots,\gamma_r}V(\gamma_1,\ldots,\gamma_r)
\sim
(T_{0}V) C_1
\sum_{\gamma_1\geq 1} u_{\gamma_1}, 
\end{eqnarray*}
which yields
\begin{eqnarray*}
\e'\log 
\sum
v_{\gamma_1,\ldots,\gamma_r}V(\gamma_1,\ldots,\gamma_r)
&=&
\e'\log 
\sum
\bar{u}_{\gamma_1,\ldots,\gamma_r}V(\gamma_1,\ldots,\gamma_r)
-\e'\log 
\sum
\bar{u}_{\gamma_1,\ldots,\gamma_r}
\\
&=&
\e'\log (T_{0}V) C_1
\sum_{\gamma_1\geq 1} u_{\gamma_1} 
- 
\e'\log
C_1
\sum_{\gamma_1\geq 1} u_{\gamma_1} 
=\log T_0 V.
\end{eqnarray*}
Taking the expectation gives (\ref{PDprop2}).
\qed

Let $\Omega$ be a set of different combinations $(i,j,l), (j,l), l$
that appear as indices of all different r.v. in Section 3.
We consider i.i.d copies $(\theta_{\omega})_{ \omega\in\Omega}$  of $\theta.$
Let $x_{\omega}(\gamma_1,\ldots,\gamma_r)$ and $v_{\gamma_1,\ldots,\gamma_r}$
be defined by (\ref{xs}) and (\ref{rweights}) and 
let $(x_{\omega})_{\omega\in\Omega}$
be an independent copy of $(x_{\omega}(1,\ldots,1))_{\omega\in\Omega}.$
The following Theorem is a consequence of
Theorem \ref{W}.

\begin{theorem}\label{rstep}
Suppose that (\ref{condition1}), (\ref{condition2}) and
(\ref{condition3}) hold and $H'_N(\vsi)$ is given by (\ref{H2}).  
Let $U_j(\eps)=U(\theta_j, x_1^j,\ldots,x_{p-1}^j, \eps),$ 
$\vec{x}=(x_1,\ldots,x_{p})$ and let $k$ be a Poisson r.v.
with mean $\e k = \alpha p.$ Then
\begin{equation}
F_N\leq \Phi_r(\zeta,m_1,\ldots,m_r)=
\log 2 + \e\log T_0\Bigl(
\Av\exp\bigl(\sum_{j\leq k} U_j(\eps) + h\eps\bigr)\Bigr)
-\alpha(p-1)\e\log T_0\la\E\ra_{\vec{x}}, 
\label{rRSB}
\end{equation}
where $T_0$ is defined in (\ref{T}) 
and $\e$ denotes the expectation w.r.t. 
$(\eta_l),$ $(\eta_l^j),$ $(\theta_j),$ $k$ and $h.$  
\end{theorem}

Of course, the Theorem implies that
$$
F_N\leq \Phi_r=\inf_{\zeta,m_1,\ldots,m_r}\Phi_r(\zeta,m_1,\ldots,m_r).
$$
One expects that as in the Sherrington-Kirkpatrick (SK) model (\cite{T3}) 
these bounds are always exact, i.e. 
$$
\lim_{N\to \infty} F_N =\inf_{r\geq 0}\Phi_r,
$$
where $\Phi_0$ was defined in (\ref{rsb00}).
Probably, this is going to be much harder to prove than
the Parisi formula in the SK model. 
The hope that these bounds are exact is not based on anything concrete 
but rather on what may be called the generalized Parisi conjecture that 
the Replica Symmetry Breaking scheme, when properly applied, 
always yields the correct free energy.

{\bf Proof of Theorem \ref{rstep}.}
The proof is almost identical to the proof of Theorem \ref{1step}
with Proposition \ref{PD2} now playing the role of Proposition \ref{PD}.
To simplify the notations we will write $\gamma=(\gamma_1,\ldots,\gamma_r)$,
and to match the notations of Theorem~\ref{W} we write 
$x_\omega^\gamma = x_\omega(\gamma_1, \cdots, \gamma_r)$.
Let us first consider the integrand
$\e\log \bigl\langle\la \E \ra_{\vec{x}^{\gamma}}\bigr\rangle$
in the last term of (\ref{Wbound}).
Let us denote 
$$
Z_{t}(\gamma)=\sum_{\vsi} \exp(-H_{N,t}^{\gamma}(\vsi))
$$
and $e(\gamma) = \la \E \ra_{\vec{x}^\gamma}.$
Note that
$e(\gamma)$ depends on $x_l(\gamma), l\geq 1$ and $Z_t(\gamma)$ depends
on $x_l^{i,j,\gamma}= x_l^{i,j}(\gamma_1, \cdots, \gamma_r), i,j,l\geq 1$ through the Hamiltonian $H_{N,t}^{\gamma}.$
By construction, given $\F_0,$
the sequences $(Z_{t}(\gamma))$ and $(e(\gamma))$ 
are defined as in (\ref{Vee}) 
and independent of each other.
Using Proposition \ref{PD2}, we get
\begin{eqnarray*}
&&
\e\log \Bigl\langle\la \E \ra_{\vec{x}^\gamma}\Bigr\rangle
=
\e \log 
\frac{\sum v_{\gamma} e(\gamma) Z_{t}(\gamma)}
{\sum v_{\gamma}Z_{t}(\gamma)}
=
\e \log 
\sum v_{\gamma} e(\gamma) Z_{t}(\gamma) -
\e \log 
\sum v_{\gamma}Z_{t}(\gamma)
\\
&&
=
\e\log  T_0(e Z_{t}) - 
\e\log  T_0 Z_{t}
=
\e\log (T_0 e) (T_0 Z_{t}) - 
\e\log T_0 Z_{t} =
\e\log T_0 e,
\end{eqnarray*}
which is independent of $t$ and, therefore,
this yields the last term of (\ref{rRSB}).

Next let us consider the first term in (\ref{Wbound}), $\varphi(0).$
First of all,
$$
\varphi(0)=\frac{1}{N}\e\log\sum_{\gamma}v_{\gamma} Z_0(\gamma),
$$
where
\begin{eqnarray*}
Z_0(\gamma)=\sum_{\vsi} \exp(-H_{N,0}^{\gamma}(\vsi))
=2^N\prod_{i=1}^{N}\Av \exp\Bigl(
\sum_{j\leq k_{i,0}} U_{i,j}^{\gamma}(\eps) + h_i\eps\Bigr)
=2^N\prod_{i=1}^{N}F_i(\gamma),
\end{eqnarray*}
and where we introduced the notation 
$F_i(\gamma)= \Av \exp\bigl(\sum_{j\leq k_{i,0}} 
U_{i,j}^{\gamma}(\eps) + h_i\eps\bigr).$ Given $\F_0,$
the sequences $(F_i(\gamma))_{\gamma\in\Natural^{r}}$
are independent for different indices $i.$
Therefore, the application of Proposition \ref{PD2} gives, writing $F_i(1)= F_i(\gamma) $ for $\gamma =(1,\cdots, 1)$, 
\begin{eqnarray*}
&&
\frac{1}{N}\e\log 2^N\sum_{\gamma}v_{\gamma}\prod_{i=1}^{N}F_i(\gamma)
=\log 2 + \frac{1}{N}\e\log T_0\prod_{i=1}^{N}F_i(1)
\\
&&
=\log 2 + \frac{1}{N}\sum_{i=1}^{N}\e\log T_0 F_i(1)
= \log 2 + \e\log T_0 F_1(1),
\end{eqnarray*}
which is precisely the first two terms in (\ref{rRSB}).

\qed

\end{document}